\documentclass[11pt]{amsart}
\usepackage{amssymb}

\usepackage{graphicx}
\usepackage{xcolor}

\usepackage{tikz}
\usetikzlibrary{snakes,arrows}

\theoremstyle{plain}
\newtheorem{theorem}{Theorem}[section]

\newtheorem{proposition}[theorem]{Proposition}

\theoremstyle{remark}

\numberwithin{equation}{section}

\begin{document}
\title[Tukey reduction]{Tukey reduction among\\ analytic directed orders
}

\author{S{\l}awomir Solecki}

\thanks{Research supported by NSF grant DMS-1001623.}

\address{Department of Mathematics\\
University of Illinois\\
1409 W. Green St.\\
Urbana, IL 61801, USA}

\email{ssolecki@math.uiuc.edu}



\begin{abstract} This is a survey of recent work on the structure of Tukey reductions
among analytic $\sigma$-ideals of compact
subsets of compact metric spaces and
analytic P-ideals of sets of natural numbers. An attempt is made to organize the results into
a unified whole. This organization makes it possible to identify natural unresolved problems.
Some new joint results with Stevo Todorcevic are announced.
\end{abstract}

\maketitle

This paper is a survey of certain results around Tukey reducibility. It is not a comprehensive survey and
I will concentrate only on what seems to be the theoretical core of the structure of Tukey reduction among definable
directed orders. Consequently, I will only touch very lightly on non-definable directed orders and on various applications of Tukey reduction.
Discussions of these topics can be found, for example, in \cite{DT}, \cite{Fr}, \cite{T1}, and the literature cited in these papers.

I will start with defining Tukey reduction. I will then describe
convenient domains for the study of Tukey reduction among definable directed orders, starting with the broadest one
and ending with six concrete representative examples. Then I will move on to describe what is known and unknown about the structure of
Tukey reductions within these domains.

The reader can find a diagram illustrating some of the material surveyed in this paper in Figure~1 in Section~\ref{S:laa}.

\section{Tukey reduction}\label{S:Tuk}

By a {\em directed order} $(D, \leq)$ we understand a partial order such that for each $x,y\in D$ there is $z\in D$
with $x,y\leq z$. Abusing notation somewhat, we will write $D$ for the directed order $(D,\leq)$. A set $A\subseteq D$ is called {\em bounded}
if there is $x\in D$ such that $y\leq x$ for each $y\in A$. The notion dual to bounded set is that of cofinal set.
A set $A\subseteq D$ is {\em cofinal} if for each $x\in D$ there $y\in A$ with $x\leq y$.

Let $D$ and $E$ be directed orders.
A function $f\colon D\to E$ is called {\em Tukey} if preimages under $f$ of sets bounded in $E$ are bounded in $D$.
We write
\[
D\leq_T E
\]
if there is a Tukey function from $D$ to $E$.
There is a notion of morphism among directed orders that is dual to Tukey morphism and
that is defined as follows.
A function $g\colon E\to D$ is {\em convergent} if images under $g$ of sets cofinal in $E$
are cofinal in $D$. It was already noted by Tukey \label{T:Tuk1} that for two directed orders $D$ and $E$,
there is a Tukey function from $D$ to $E$ if and only if there is a convergent function from $E$ to $D$.

If two directed orders are Tukey reducible to each other, we say that they are {\em Tukey equivalent}.
We write
\[
D\equiv_T E
\]
if both $D\leq_TE$ and $E\leq_TD$. As shown by Tukey \cite{Tu}, this condition
can be phrased in a way that does not involve Tukey reduction.

\begin{theorem}[\cite{Tu}]
Let $D$ and $E$ be directed orders.
Then $D\equiv_T E$ if and only if there is a directed order $F$
such that $D$ and $E$ embed into $F$ as cofinal subsets.
\end{theorem}

Tukey reduction was originally introduced \cite{Tu} in the theory of net convergence in general topological spaces
in order to formulate the important notion of subnet.
Later Isbell~\cite{Is} realized that Tukey reduction can be fruitfully used to compare
directed orders coming from topology and analysis.
Building on this insight and on the insight of Schmidt~\cite{Sc}, who connected Tukey reduction with the cardinals of additivity and cofinality
of a directed order, Fremlin~\cite{Fr} employed Tukey reduction to explain inequalities between certain cardinal invariants of
the continuum. This explanation is based on the following observation. Let $D$ be a directed order. Define additivity of $D$
\[
{\rm add}(D)
\]
to be the minimal cardinality of an unbounded subset of $D$, and let cofinality of $D$
\[
{\rm cof}(D)
\]
be the minimal cardinality of a cofinal subset of $D$. Now, it is an observation going back to Schmidt~\cite{Sc}  that if $D\leq_T E$, then
\[
{\rm add}(E)\leq {\rm add}(D)\hbox{ and }{\rm cof}(D)\leq {\rm cof}(E).
\]

In applications to cardinal invariants of the continuum, a notion slightly weaker than Tukey reduction is also of interest. Let $D$ and $E$ be directed orders. We
write $D\leq_T^\omega E$ if there is a function $f\colon D\to E$ such that preimages of $\sigma$-bounded sets are
$\sigma$-bounded. (A set is {\em $\sigma$-bounded} if it is a countable union of bounded sets.) It is clear that $D\leq_TE$ implies $D\leq_T^\omega E$.
We define
\[
{\rm add}^\omega(D)
\]
as the smallest cardinality of a non-$\sigma$-bounded subset of $D$. It is now easy to see \cite{Fr} that
$D\leq_T^\omega E$ implies that
\[
{\rm add}^\omega(E)\leq {\rm add}^\omega(D)\;\hbox{ and }\;{\rm cof}(D)\leq \max(\omega, {\rm cof}(E)).
\]
We write
\[
D\equiv_T^\omega E
\]
if both $D\leq_T^\omega E$ and $E\leq_T^\omega D$.

Using Tukey reduction as a means of comparison among directed orders occurring in topology or analysis has
some interesting advantages. On the one hand, Tukey reduction is an abstract and general notion, so it makes it possible to compare very
diverse directed orders and, even in situations when it is applied to compare directed orders related to each other, it puts such
comparisons in a broad context. For example,
one class of directed orders studied quite extensively consists
of ideals of sets of natural numbers taken with inclusion as the directed order relation.
There are several useful notions that
allow comparison between such ideals: the Rudin--Keisler, Rudin--Blass, or Kat{\v e}tov reductions. All of them, however, depend on the fact that the ideals
are defined on the set of natural numbers. Tukey reduction provides a ``coordinate free" way of comparing such ideals and places such comparisons against a
larger backdrop.
On the other hand, remarkably, despite Tukey reduction being an abstract notion,
in controlled situations, the existence of an abstract Tukey reduction implies the existence of a {\em definable} such reduction; see Theorem~\ref{T:deft}.

\section{Basic orders, ideals, six examples}\label{S:ide}

\subsection{Basic orders} The theory of Tukey reduction can be nicely developed in a class of directed orders whose
underling sets are appropriately topologized. This class contains the main examples of definable directed orders. The following definition
is due to Solecki and Todorcevic and comes from \cite{ST1}.
A directed order $D$ is called {\em basic} if
\begin{enumerate}
\item[---] $D$ is a separable metric space;

\item[---] each two elements of $D$ have
the least upper bound and the operation of taking the least upper bound is a continuous function from $D\times D$ to $D$;

\item[---] each bounded
sequence has a convergent subsequence;

\item[---] each convergent sequence has a bounded subsequence.
\end{enumerate}

It was pointed out by Fremlin \cite[Proposition 513K]{Fr:mea} that the topology on a basic order is determined by
the order relation.

Basic orders whose underlying topology is analytic will be called {\em analytic basic orders}. (A metric separable space is {\em analytic} if
it is a continuous image of a Polish space.)
There are two ``self-improvement" results for analytic basic orders. The first one of which concerns the topology on a basic order.

\begin{theorem}[\cite{ST1}]\label{T:pol}
Let $D$ be a basic order. If the topology on $D$ is analytic, then it is Polish.
\end{theorem}

\noindent The above theorem was anticipated in the results of Christensen~\cite[Theorem 3.3]{Ch} and of
Kechris--Louveau--Woodin~\cite{KLW}, who proved that each analytic $\sigma$-ideal of compact subsets
of a compact metric space is $G_\delta$.

The second of the ``self-improvement"
results concerns morphisms and shows that Tukey reducibility among analytic basic orders can be always
witnessed by definable functions.

\begin{theorem}[\cite{ST1}]\label{T:deft}
Let $D$ and $E$ be analytic basic orders. Assume $D\leq_T E$. Then
there exist a Tukey function from $D$ to $E$
that is measurable with respect to the $\sigma$-algebra generated by analytic sets.
\end{theorem}

Below we will be interested exclusively in analytic basic orders. (Note, however, that
non-analytic basic orders have been investigated in the literature; see for example \cite{DT}.)
This class is described by imposing
a definability condition but, as stated in the next theorem, it forms an initial set of basic orders; so the definable limitation turns out
to be a complexity limitation as well.

\begin{theorem}[\cite{ST1}]
Let $D$ and $E$ be basic orders with $D\leq_TE$. If $E$ is analytic, then so is $D$.
\end{theorem}

To concentrate on non-trivial basic orders, we will be interested in analytic basic orders that are
not locally compact. This class has a minimal element.
Anticipating Subsection~\ref{Su:ex}, we introduce here an analytic non-locally compact basic order.
Consider the set ${\mathbb N}^{\mathbb N}$ of all functions from $\mathbb N$ to $\mathbb N$ with the product topology and the pointwise inequality between functions.

\begin{theorem}[\cite{ST1}]\label{T:nn}
Let $D$ be an analytic non-locally compact basic order. Then ${\mathbb N}^{\mathbb N}\leq_T D$.
\end{theorem}

\noindent
The special case of the above theorem for analytic P-ideals was proved earlier by Todorcevic~\cite{T2}.

\subsection{Ideals}
We shift our attention to two subclasses of analytic basic orders from which most natural examples come.
Each of these two classes consist of ideals taken with inclusion as directed order. The reader may consult two
recent surveys
\cite{Hr} and \cite{MZ} for more background information on ideals.

First, we note, somewhat
academically, that there is no loss of generality in considering only ideals rather than general directed orders,
as each directed order is easily seen to
be Tukey equivalent to the ideal of its bounded subsets.
However, we will not consider ideals in full generality but rather limit our attention to ideals that are also analytic
basic orders.
As is often the case in mathematics, this domain of investigation
splits into the compact and the discrete subdomains. More precisely, we will be interested in analytic $\sigma$-ideals of compact sets and
analytic P-ideals of subsets of $\mathbb N$, both taken with inclusion as partial order and with appropriate topologies. We describe the two classes in turn.

Consider a compact metric space $X$ and equip the space ${\mathcal K}(X)$ of all compact subsets of $X$ with the usual
Vietoris topology, which makes ${\mathcal K}(X)$ into a compact metric space. A set ${\mathcal I}\subseteq {\mathcal K}(X)$ is a
{\em $\sigma$-ideal of compact sets} if it is closed under taking compact subsets and countable compact unions. It is easy to see
that a $\sigma$-ideal of compact sets with inclusion and the topology inherited from ${\mathcal K}(X)$ is a basic order. It is
an analytic basic order if the Vietoris topology on it is analytic. It was proved by Kechris--Louveau--Woodin~\cite{KLW}, and follows also from
Theorem~\ref{T:pol}, that in that case
$\mathcal I$ is a $G_\delta$ subset of ${\mathcal K}(X)$. In certain situations, a somewhat more general notion of a {\em relative $\sigma$-ideal
of compact sets} is natural and useful; see \cite{ST1} or \cite{Ma3}. However, we will not consider this generalization below.

To describe the other class of ideals, consider the powerset ${\mathcal P}({\mathbb N})$ of $\mathbb N$ identified with $2^{\mathbb N} = \{ 0,1\}^{\mathbb N}$
and through this identification  equipped with the usual compact product
topology. A set $I\subseteq {\mathcal P}({\mathbb N})$ is an ideal if it is closed under taking finite unions and subsets. It can be checked
that $I$ taken with inclusion and with the topology inherited from ${\mathcal P}({\mathbb N})$ is an analytic basic order
precisely when $I = {\mathcal P}(x)$ for some $x\subseteq {\mathbb N}$, so it would appear that this class  contains only trivial
examples. As it turns out, this is because we were too simple-minded in our choice of topology.
One does get a large number of important examples of ideals of subsets of
$\mathbb N$ that are analytic basic orders
as follows. An ideal $I$ is called a {\em P-ideal} if for each sequence $x_n\in I$, $n\in {\mathbb N}$, there is $x\in I$ such that
$x_n\setminus x$ is finite for each $n$. It was proved in \cite{So1} that an ideal $I$ is a P-ideal that is analytic with the topology inherited
from $2^{\mathbb N}$ if and only if there exists a lower semicontinuous submeasure $\phi\colon {\mathcal P}({\mathbb N})\to [0,\infty]$ such that
\[
I = {\rm Exh}(\phi)= \{ x\in {\mathcal P}({\mathbb N})\colon \lim_n\phi(x\setminus \{ 0, 1,\dots, n\})=0\}.
\]
We can always assume that $0< \phi(\{ n\})<\infty$ for each $n\in {\mathbb N}$.
It is easy to see that such an $I$ becomes a Polish space with the submeasure topology given by the metric
\[
d_\phi(x,y) = \phi((y\setminus x)\cup (x\setminus y)).
\]
It is easy to see that if ${\rm Exh}(\phi) = {\rm Exh}(\phi')$  for two lower semicontinuous
submeasures $\phi$ and $\phi'$, then
for each $\epsilon>0$ there is $\epsilon'>0$ such that $\phi'(x)<\epsilon'$ implies $\phi(x)<\epsilon$ for each subset $x$ of $\mathbb N$
and vice versa. Consequently, the submeasure topology does
not depend on a particular choice of $\phi$; it is determined by $I$.
It is easy to check that $I$ with inclusion and with this topology is an analytic basic order.

This definition may strike the reader not acquainted with the area as not entirely natural. However, first, many important examples
are of this form and, second, this is the only way of making ideals of subsets of $\mathbb N$ into analytic basic orders. Indeed,
it was proved in \cite{ST1} that if an ideal $I\subseteq {\mathcal P}({\mathbb N})$, taken with inclusion and with a topology $\tau$ containing
the topology inherited from ${\mathcal P}({\mathbb N})$, is an analytic basic order, then $I$ is an analytic P-ideal and $\tau$
is the submeasure topology.

\subsection{Examples}\label{Su:ex} We move now to a description of some concrete examples of ideals. They come from
different important classes of
ideals. We omit many interesting other examples, which can be easily found in the literature cited below.

The first four examples were already considered by Isbell~\cite{Is}.

1. The basic order ${\mathbb N}^{\mathbb N}$ was introduced when stating Theorem~\ref{T:nn}.
This partial order can be considered as both an analytic $\sigma$-ideal of compact subsets of a compact metric space
and an analytic P-ideal of subsets of $\mathbb N$. So ${\mathbb N}^{\mathbb N}$ is Tukey equivalent with the $\sigma$-ideal
of all compact subsets of
$[0,1]$ not containing rational numbers and it is also Tukey equivalent with the P-ideal consisting of subsets of
${\mathbb N}\times {\mathbb N}$ contained in the subgraph of a function from
$\mathbb N$ to $\mathbb N$.

2. Let
\[
{\mathcal Z}_0 = \{ x\in {\mathcal P}({\mathbb N}) \colon \lim_n \frac{|x\cap \{ 0, 1, \dots, n\}|}{n+1} =0\}.
\]
This is the P-ideal consisting of density zero subsets of $\mathbb N$. It is easy to see that the following
lower semicontinuous submeasure
\[
\phi_0(x) = \sup_n \frac{|x\cap \{ 0, 1, \dots, n\}|}{n+1}
\]
is such that ${\mathcal Z}_0= {\rm Exh}(\phi_0)$.

3. Let
\[
\ell_1 = \{ x\in {\mathcal P}({\mathbb N})\colon \sum_{n\in x}\frac{1}{n+1}<\infty\}.
\]
Again it is easy to check that $\ell_1$ is a P-ideal. Moreover, one checks that the lower semicontinuous
submeasure
\[
\phi_1(x) =\sum_{n\in x} \frac{1}{n+1}
\]
gives $\ell_1 = {\rm Exh}(\phi_1)$. The ideal $\ell_1$ is Tukey equivalent with the directed order of all summable
sequences of non-negative real numbers taken with pointwise inequality as the order relation \cite{Fr}. This
observation explains the symbol used to denote the ideal.

4. Let
\[
{\rm NWD}
\]
be the ideal of all compact nowhere dense subsets of $2^{\mathbb N}$. A straightforward argument shows that
the ideal $\rm NWD$ is an analytic $\sigma$-ideal.

5. A new type of an analytic $\sigma$-ideal of compact subsets of $2^{\mathbb N}$ was discovered recently by M{\'a}trai \cite{Ma1}. We present
a $\sigma$-ideal of this type below. This particular $\sigma$-ideal comes from \cite{So2}
and is easier to describe than the original $\sigma$-ideal from \cite{Ma1}.

We consider sequences ${\bar s} = (s_0, s_1, \dots)$ that are infinite or finite with an even number of entries, where each $s_i$
is a function from a non-empty finite subset of $\mathbb N$ to $2=\{ 0,1\}$, and where for each $i$ each element of the domain of $s_i$ is less
than each element of the domain of $s_{i+1}$. Let $\mathcal R$ be the set of all such sequences.
For ${\bar s}\in {\mathcal R}$, define
\[
[{\bar s}] = \{ x\in 2^{\mathbb N}\colon s_{2i}\subseteq x\hbox{ or }s_{2i+1}\subseteq x\hbox{ for each }i\}.
\]
Define
\[
{\mathcal I}_0 = \{ K\in {\mathcal K}(2^{\mathbb N})\colon K\cap [{\bar s}] \hbox{ is nowhere dense in }
[{\bar s}] \hbox{ for each } {\bar s}\in {\mathcal R}\; \}.
\]
One checks that ${\mathcal I}_0$ is an analytic $\sigma$-ideal \cite{So2}.

The examples listed above represent all the classes of ideals relevant in the sequel: $\rm NWD$ is an analytic $\sigma$-ideal of compact sets with
property $(*)$, ${\mathcal I}_0$ is an analytic $\sigma$-ideal of compact sets without $(*)$, ${\mathcal Z}_0$ is a density-like analytic P-ideal of subsets
of $\mathbb N$, $\ell_1$ is an analytic P-ideal of subsets of $\mathbb N$ that is not density-like.
(Property $(*)$ and density-like will be defined later.)
Additionally, ${\mathbb N}^{\mathbb N}$ is the unique up
to Tukey equivalence, non-locally compact basic order that is both an analytic $\sigma$-ideal of compact sets and an analytic P-ideal of subsets of $\mathbb N$;
see Theorem~\ref{T:psig}.

There exist, however, many natural basic orders that were omitted from the list above.
One of them has been carefully studied and for this reason we will have it here as our last example.

6. Let
\[
{\mathcal E}_\mu
\]
be the analytic $\sigma$-ideal of all compact Lebesgue measure zero subsets of the unit interval taken with inclusion.

We should point out that other natural partial orders that are not themselves P-ideals or $\sigma$-ideals or even
basic orders can be analyzed in the set-up described above. For example, one of those is the directed order whose underlying
set is the family $\rm MGR$ of meager subsets of $[0,1]$ and whose order is inclusion. As shown in \cite{Fr} we have
\[
{\rm MGR} \equiv_T^\omega {\rm NWD}.
\]
And, by the remarks in Section~\ref{S:Tuk}, the study of cardinal invariants of additivity and cofinality of meager sets boils down to the study of
cardinal invariants of $\rm NWD$. Similarly, consider the set $\rm NULL$ of all Lebesgue null subsets of the interval $[0,1]$ taken with inclusion.
Building on earlier work of Bartoszy{\'n}ski~\cite{Ba} and Raisonnier--Stern~\cite{RS}, Fremlin showed in \cite{Fr}, that
\[
{\rm NULL}\equiv_T^\omega \ell_1.
\]
And again the study of additivity and cofinality of $\rm NULL$ is reduced to the study of these cardinal invariants of the
analytic basic order $\ell_1$.

\section{Structure of the classes of $\sigma$-ideals and P-ideals}

Analytic $\sigma$-ideals of compact sets and analytic $P$-ideals of subsets of $\mathbb N$ are locally compact with the topologies
that make them basic orders only in trivial situations. Indeed, it was proved in \cite{KLW} that a $\sigma$-ideal of compact subsets
of a compact metric space $X$ is locally compact with the Vietoris topology precisely
when it is the family of all compact subsets of $U$ for a fixed open subset $U$ of $X$.
It was proved in \cite{So1} that an analytic P-ideal is locally compact with its submeasure topology precisely when it is of the form
$\{ x\subseteq {\mathbb N}\colon x\cap a\hbox{ is finite}\}$ for some fixed $a\subseteq {\mathbb N}$. In both situations the resulting
directed orders are Tukey equivalent to the one element order or to $\mathbb N$ taken with the usual inequality relation.
{\em To avoid these trivial situations,
from this point on, we consider only
non-locally compact analytic $\sigma$-ideals of
compact subsets of a compact metric space and
non-locally compact analytic P-ideals of subsets of $\mathbb N$.
For simplicity, we refer to the former as {\bf $\sigma$-ideals} and to the latter {\bf P-ideals}.}

The following theorem due to Solecki and Todorcevic~\cite{ST1} shows that there are essentially no Tukey reductions from the P-ideal side to
the $\sigma$-ideal side. On the other hand, we will see later that there do exist Tukey reductions in the opposite direction.

\begin{theorem}[\cite{ST1}]\label{T:psig}
Let $I$ be a P-ideal and let $\mathcal I$ be a $\sigma$-ideal. Then $I\not\leq_T {\mathcal I}$ unless $I$ is isomorphic
to the P-ideal of subsets of ${\mathbb N}\times {\mathbb N}$ contained in the subgraph of a function from ${\mathbb N}$ to ${\mathbb N}$, so
$I\equiv_T {\mathbb N}^{\mathbb N}$.
\end{theorem}

\noindent A particular instance of the above theorem, namely ${\mathcal Z}_0\not\leq_T {\rm NWD}$, was proved earlier by Fremlin in \cite{Fr}.

Louveau and Veli{\v c}kovi{\'c}~\cite{LV} and, independently, Todorcevic~\cite{T} noticed that the following theorem is a consequence of a general result
of Fremlin~\cite{Fr} on producing Tukey reduction to $\ell_1$ and the representation of P-ideals as ${\rm Exh}(\phi)$ from \cite{So1}.

\begin{theorem}[\cite{LV}, \cite{T}]
$\ell_1$ is largest with respect to Tukey reduction among all P-ideals, that is,
$I\leq_T \ell_1$ for each P-deal $I$.
\end{theorem}

It is not known if the class of $\sigma$-ideals has a largest element with respect to Tukey reduction.

There are important subclasses of $\sigma$-ideals and P-ideals that exhibit a higher degree of additivity.
They are in some sense analogous to each other. We describe them below.

On the side of $\sigma$-ideals we find the following notion introduced in \cite{So2}. Let $\mathcal I$ be a $\sigma$-ideal. Let $X$
be the compact metric space underlying $\mathcal I$.
We say that $\mathcal I$ has {\em property $(*)$} if for each sequence $(K_n)$ of sets in $\mathcal I$ there
is a $G_\delta$ subset $G$ of $X$ such that $\bigcup_n K_n \subseteq G$ and all compact subsets of $G$ are in $\mathcal I$.
One checks that $\rm NWD$ and ${\mathcal E}_\mu$ have property $(*)$ and that ${\mathcal I}_0$ does not; see \cite{Ma1} and \cite{So2}.
As argued in \cite{So2} all ``naturally occurring" $\sigma$-ideals have property $(*)$. In fact, the $\sigma$-ideal constructed
in \cite{Ma1} was the first example of a $\sigma$-ideal without $(*)$.

On the P-ideal side, the following
definition was introduced by Solecki and Todorcevic in \cite{ST1}. Let $I={\rm Exh}(\phi)$ be a P-ideal for
a lower semicontinuous submeasure $\phi$. We call $I$ {\em density-like} if for each $\epsilon>0$ there is $\delta>0$ such that
for each sequence $(x_n)$ of sets in $I$ with $\phi(x_n)\leq \delta$ there is an infinite set $b\subseteq {\mathbb N}$ with
\[
\phi(\bigcup_{n\in b} x_{n})\leq \epsilon.
\]
Whether or not $I$ is density-like depends only on $I$ and not on
the choice of $\phi$ with $I={\rm Exh}(\phi)$.
One checks without difficulty that ${\mathcal Z}_0$ is density-like and that $\ell_1$ is not.

It was proved by Solecki~\cite{So2} that among $\sigma$-ideals with $(*)$ there is a top element.

\begin{theorem}[\cite{So2}]
$\rm NWD$ is largest with respect to Tukey reduction among all $\sigma$-ideals with property $(*)$,
that is, ${\mathcal I}\leq_T {\rm NWD}$ for each $\sigma$-ideal with $(*)$ $\mathcal I$.
\end{theorem}

It is not know if there is a largest density-like P-ideal.

The following theorem, due to Louveau and Veli{\v c}kovi{\'c}~\cite{LV} and M{\'a}trai~\cite{Ma3},
illustrates richness of Tukey reduction among P-ideals.

\begin{theorem}
\begin{enumerate}
\item[(i)] {\rm (\cite{LV})}
The partial order ${\mathcal P}({\mathbb N})/{\rm Fin}$ with almost inclusion embeds into the class of density-like P-ideals with $\leq_T$.

\item[(ii)] {\rm (\cite{Ma3})}
The partial order ${\mathcal P}({\mathbb N})/{\rm Fin}$ with almost inclusion embeds into the class of non-density-like P-ideals with $\leq_T$.
\end{enumerate}
\end{theorem}

\noindent Strictly speaking, in \cite{Ma3}, it is shown that the class of P-ideals that are $F_\sigma$ subsets of
 $2^{\mathbb N}$ is rich in the way described in point (ii), but it is proved in \cite{ST1} that P-ideals that are
 $F_\sigma$ are not density-like.

The corresponding result for $\sigma$-ideals are not known, though some progress on this question was made in \cite{Ma3}. In fact,
one of the most interesting problems in the area may be the challenge of sorting out the structure of Tukey reduction among the many mathematically
natural $\sigma$-ideals with property $(*)$. A list of such ideals can be found in \cite[Section~2]{So2}.

\section{Tukey reductions among the examples}\label{S:laa}

We now go back to the concrete examples from Section~\ref{S:ide}: ${\mathbb N}^{\mathbb N}$, ${\mathcal Z}_0$, $\ell_1$, ${\rm NWD}$, ${\mathcal I}_0$,
and ${\mathcal E}_\mu$.
The structure of Tukey reduction among these ideals has been completely determined.
The first theorem shows the Tukey comparisons between elements of the same class: P-ideals and $\sigma$-ideals.

\begin{theorem}\label{T:brrr}
\begin{enumerate}
\item[(i)] {\rm (\cite{Fr}, \cite{Ma2},  \cite{Ma3}, \cite{MS}, \cite{ST2})}
\[
{\mathbb N}^{\mathbb N} <_T {\mathcal E}_\mu<_T{\rm NWD} <_T {\mathcal I}_0.
\]

\item[(ii)] {\rm (\cite{Fr}, \cite{Is}, \cite{LV})}
\[
{\mathbb N}^{\mathbb N} <_T {\mathcal Z}_0 <_T \ell_1
\]
\end{enumerate}
\end{theorem}

In point (i), ${\mathbb N}^{\mathbb N} <_T {\mathcal E}_\mu$ and ${\mathcal E}_\mu\leq_T {\rm NWD}$ were
proved by Fremlin~\cite{Fr} with strictness of the latter inequality following
from Theorem~\ref{T:acr}(ii) and (iii) below and established in \cite{Fr}, \cite{Ma3}, and \cite{ST2}. The inequality ${\rm NWD} <_T {\mathcal I}_0$ in point (i)
was proved independently by M{\'a}trai \cite{Ma2} and by  Moore and Solecki~\cite{MS}.
In point (ii), ${\mathbb N}^{\mathbb N} <_T {\mathcal Z}_0$ is due to Isbell~\cite{Is},  ${\mathcal Z}_0 \leq_T \ell_1$ to Fremlin~\cite{Fr},
and $\ell_1\not\leq_T {\mathcal Z}_0$
to Louveau and Veli{\v c}kovi{\' c}~\cite{LV}.
In point (ii), strictness of the inequality ${\mathcal Z}_0<_T \ell_1$ also
follows from the general Theorems~\ref{T:htin} and
\ref{T:htpr} as ${\mathcal Z}_0$ is density-like.

Now we take a look at comparisons between elements of different classes.
It follows from the general Theorem~\ref{T:psig} that
${\mathcal Z}_0$ and $\ell_1$ are not Tukey reducible to ${\rm NWD},\, {\mathcal I}_0,$ or ${\mathcal E}_\mu$.
In the opposite direction, we have the following theorem.

\begin{theorem}\label{T:acr}
\begin{enumerate}
\item[(i)] {\rm (\cite{Fr})} ${\rm NWD}<_T \ell_1$

\item[(ii)] {\rm (\cite{Fr})} ${\mathcal E}_\mu<_T {\mathcal Z}_0$;

\item[(iii)] {\rm (\cite{Ma3}, \cite{ST2})} ${\rm NWD}\not\leq_T {\mathcal Z}_0$

\item[(iv)] {\rm (\cite{Ma2})} ${\mathcal I}_0\not\leq_T\ell_1$
\end{enumerate}
\end{theorem}

In the above theorem, points (i) and (ii) are due to Fremlin~\cite{Fr}.
Point (iii) was proved independently by M{\'a}trai~\cite{Ma3} and Solecki--Todorcevic~\cite{ST2}. Point (iv) was proved by M{\'a}trai in \cite{Ma2}. Strictly
speaking, it was not proved there for
${\mathcal I}_0$, but for the $\sigma$-ideal defined in \cite{Ma1}. This proof can be adapted to yield point (iv) for
the $\sigma$-ideal ${\mathcal I}_0$ described here in Subsection~\ref{Su:ex}. Note that strictness of the inequalities in points (i) and (ii) follows
also from the general Theorem~\ref{T:psig}.

Theorems~\ref{T:brrr} and \ref{T:acr} together determine all inequalities and lack thereof among the examples. Note also that by the discussion
earlier on Theorem~\ref{T:acr}(i) implies that
\begin{equation}\label{E:numg}
{\rm add(NULL)} \leq {\rm add(MGR)}\hbox{ and }{\rm cof(MGR)}\leq {\rm cof(NULL)}.
\end{equation}

The diagram below shows the basic orders considered in this paper and their placement within subclasses of basic orders. 
Tukey reductions among them exist precisely when indicated by arrows or compositions of arrows.

\medskip

\begin{center}
\begin{tikzpicture}
    \tikzset{>=latex'}
    \draw (0,0) -- (-4,3) -- (-4,6) -- (-1.5,6) -- (-1.5,3) -- cycle;
    \draw (0,0) -- ( 4,3) -- ( 4,6) -- ( 1.5,6) -- ( 1.5,3) -- cycle;
    \draw (-4,3) -- (-1.5,3);
    \draw ( 4,3) -- ( 1.5,3);
    \filldraw (0,0) circle (1.5pt) node[anchor=north] {\( \mathbb{N}^{\mathbb{N}} \)};
    \filldraw (-2.75,3) circle (1.5pt) node[anchor=north] {\( \mathrm{NWD} \)};
    \filldraw (-2.75,5) circle (1.5pt) node[anchor=south] {\( \mathcal{I}_{0} \)};
    \filldraw (-1.8,2) circle (1.5pt) node[anchor=east] {\( \mathcal{E}_{\mu} \)};
    \filldraw ( 2.2,2.5) circle (1.5pt) node[anchor=west] {\( \mathcal{Z}_{0} \)};
    \filldraw ( 2.75,6) circle (1.5pt) node[anchor=north] {\( \ell_{1} \)};
    \draw[thin,->] (0,0) -- (-1.7,2*17/18);
    \draw[thin,->] (-1.8,2) -- (-1.8-0.6*0.95,2+0.6*1);
    \draw[thin,->] (0,0) -- ( 0.94*2.2,0.94*2.5);
    \draw[thin,->] (2.2,2.5) -- ( 2.2+0.85*0.55,2.5+0.85*3.5);
    \draw[thin,->] (-2.75,3) -- (-2.75,4.9);
    \draw[thin,->] (-2.75,3) -- (-2.75 + 0.94*5.5,3 + 0.94*3);
    \draw[thin,->] (-1.8,2) -- (-1.8 + 0.95*4,2 + 0.95*0.5);
    \draw[decoration={
      brace,
    },
    decorate] (-0.7,0) -- (-4,3-3*0.7/4) node[pos=0.5, anchor=north, yshift=-1mm,rotate=-36.87]
    {property \( (*) \)};
    \draw[decoration={
      brace,
      mirror,
    },
    decorate] (0.7,0) -- (4,3-3*0.7/4) node[pos=0.5, anchor=north, yshift=-1mm,rotate=36.87]
    {density-like};
    \draw (-2.75,6.3) node {\( \sigma \)-ideals};
    \draw ( 2.75,6.3) node {\( P \)-ideals};
\end{tikzpicture}

\smallskip

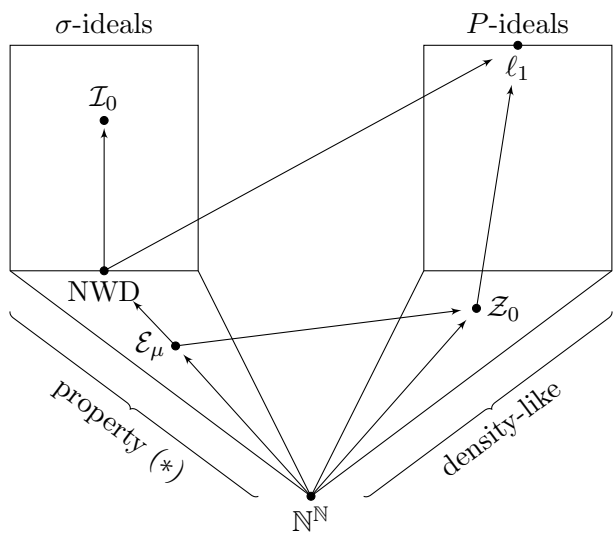
\centerline{\textsc{Figure 1.} Basic orders}
\end{center}

\medskip


Partly because of application \eqref{E:numg} to cardinal invariants, Theo\-rem~\ref{T:acr}(i) appears to be the most intriguing Tukey reduction among
the ideals considered here. It immediately leads to the problem of
characterizing those P-ideals $I$ for which ${\rm NWD}\leq_TI$. The first step in this direction is the following theorem due to
Solecki and Todorcevic~\cite{ST2}  that extends Theorem~\ref{T:acr}(iii).

\begin{theorem}[\cite{ST2}]\label{T:dens}
Let $I$ be a density-like P-ideal. Then ${\rm NWD}\not\leq_T I$.
\end{theorem}

We report below some results related to the problem of characterizing those P-ideals $I$ for which $I\leq_T {\rm NWD}$. These results were
obtained by Todorcevic and the author.

We define a new rank with values in $\omega_1+1$ on P-ideals. Note that a rank with values in $\omega_1+1$ on
$\sigma$-ideals was described in \cite{So2}.
Let a P-ideal $I$ be represented as $I= {\rm Exh}(\phi)$ for a lower semicontinuous submeasure $\phi$.
Given a sequence $(x_n)$ of sets in $I$ and $\epsilon >0$, the set
\[
K_\epsilon=\{ b\subseteq {\mathbb N}\colon \phi(\bigcup_{n\in b} x_n) \leq \epsilon \} \subseteq 2^{\mathbb N}
\]
is compact. There are two possibilities for the iteration of the Cantor--Bendixson derivative applied to this set.
Either there is a smallest countable ordinal $\alpha$ such that the $\alpha+1$ Cantor--Bendixson derivative of
$K_\epsilon$ is empty or all Cantor--Bendixson derivatives of $K_\epsilon$ are non-empty. In the first case, let
${\rm height}(K_\epsilon)$ be equal to $\alpha$; in the second case, let it be equal to  $\omega_1$.
Let now $\epsilon,\, \delta>0$ and $\alpha < \omega_1$ be given. We say that
{\em $P_{\epsilon,\delta}(\alpha)$ holds for $I$} if for every sequence $(x_n)$ of sets in $I$ with $\phi(x_n)\leq\delta$
\[
{\rm height}(\{ b \subseteq {\mathbb N}\colon \phi(\bigcup_{n\in b} x_n) \leq \epsilon \}) \geq \alpha.
\]
Define
\[
{\rm ht}(I) = \min \{\alpha\in\omega_1\colon \exists \epsilon>0\, \forall \delta>0\; P_{\epsilon,\delta}(\alpha) \hbox{ fails}\}
\]
if the set on the right hand side is non-empty, and let
\[
{\rm ht}(I) = \omega_1
\]
otherwise.

\begin{proposition}[Solecki--Todorcevic]
Let $I$ be a P-ideal. Then
\begin{enumerate}
\item[(i)] ${\rm ht}(I)$ does not depend on the choice of a lower semicontinuous submeasure $\phi$
with $I= {\rm Exh}(\phi)$;

\item[(ii)]
\[
{\rm ht}(I) = \omega^{\omega^\alpha},\hbox{ for some }\alpha<\omega_1,\hbox{ or }\; {\rm ht}(I)=\omega_1.
\]
\end{enumerate}
\end{proposition}

The following result, essentially proved in \cite{ST2} (where, however, the notion of ${\rm ht}$ was missing),
gives a characterization of P-ideals with  the largest and
smallest values of height. The proof of point (i) uses ideas of Fremlin~\cite{Fr}.

\begin{theorem}[Solecki--Todorcevic]~\label{T:htpr}
Let $I$ be a P-ideal. Then
\begin{enumerate}
\item[(i)] ${\rm ht}(I) =\omega$ if and only if $I\equiv_T \ell_1$;

\item[(ii)] ${\rm ht}(I)=\omega_1$ if and only if $I$ is density-like.
\end{enumerate}
\end{theorem}

The next result shows that height is an invariant of Tukey reduction.

\begin{theorem}[Solecki--Todorcevic]~\label{T:htin}
Let $I,J$ be P-ideals. If $I\leq_T J$, then ${\rm ht}(J)\leq {\rm ht}(I)$.
\end{theorem}

It follows from Theorem~\ref{T:htin} and Theorem~\ref{T:htpr}(ii) that the class of density-like P-ideals
is closed downward, a result established already in \cite{ST2}.

It is clear from what was said above that there is a minimal ordinal $\beta$ such that ${\rm NWD}\leq_T I$ for all P-ideals $I$
with ${\rm ht}(I)<\beta$ and that we have $\omega^\omega\leq\beta\leq \omega_1$ for this ordinal. It seems likely that
$\beta$ is equal to one of the two extreme values.

\medskip

\noindent {\bf Acknowledgement.} I would like to thank Stevo Todorcevic for his comments on an earlier version of the paper and
Jordi Lopez Abad and Kostya Slutsky for preparing the drawing for me.


\begin{thebibliography}{10}
\bibitem{Ba} T. Bartoszy{\'n}ski, {\em Additivity of measure implies additivity of category}, Trans. Amer.
Math. Soc. 281 (1984), 209--213.

\bibitem{Ch} J.P.R. Christensen, {\em Topology and Borel Structure}, North-Holland Mathematics Studies, 10, Notas de Matem‡tica, 51,
North-Holland, American Elsevier, 1974.

\bibitem{DT} N. Dobrinen, S. Todorcevic, {\em Tukey types of ultrafilters}, Illinois J. Math. to appear.

\bibitem{Fr} D.H. Fremlin, {\em The partially ordered sets of measure theory and Tukey's ordering}, Note
di Matematica 11 (1991), 177--214.

\bibitem{Fr:mea} D.H. Fremlin, {\em Measure Theory}, Vol. 5, Torres Fremlin, 2008.

\bibitem{Hr} M. Hru{\v s}{\' a}k,
{\em Combinatorics of filters and ideals}, in {\em Set Theory and Its Applications}, 29--69,
Contemp. Math., 533, Amer. Math. Soc., 2011.

\bibitem{Is} J.R. Isbell, {\em Seven cofnal types}, J. London Math. Soc. 4 (1972), 651--654.

\bibitem{KLW} A.S. Kechris, A. Louveau, W.H. Woodin,
{\em The structure of ?-ideals of compact sets}, Trans. Amer. Math. Soc. 301 (1987), 263--288.

\bibitem{LV} A. Louveau, B. Veli{\v c}kovi{\'c}, {\em Analytic ideals and cofnal types}, Ann. Pure Appl. Logic
99 (1999), 171--195.

\bibitem{MZ} {\'E}. Matheron, M. Zelen{\'y}, {\em Descriptive set theory of families of small sets}, Bull. Symbolic Logic 13 (2007), 482--537.

\bibitem{Ma1} T. M{\'a}trai, {\em Kenilworth}, Proc. Amer. Math. Soc. 137 (2009), 1115--1125.

\bibitem{Ma2} T. M{\'a}trai, {\em On a $\sigma$-ideal of compact sets}, Topology Appl. 157 (2010), 1479--1484.

\bibitem{Ma3} T. M{\'a}trai, {\em More cofinal types of definable directed orders}, to appear.

\bibitem{MS} J.T. Moore, S. Solecki, {\em A $G_\delta$ ideal of compact sets strictly above the nowhere dense
ideal in the Tukey order}, Ann. Pure Appl. Logic 156 (2008), 270--273.

\bibitem{RS} J. Raissonier, J. Stern, {\em The strength of measurability hypotheses}, Israel J. Math. 50
(1985), 337--349.

\bibitem{Sc} J. Schmidt, {\em Konfinalit{\" a}t},  Z. Math. Logik Grundlagen Math. 1 (1955), 271--303.

\bibitem{So1} S. Solecki, {\em Analytic ideals and their applications}, Ann. Pure Appl. Logic 99 (1999), 51--72.

\bibitem{So2} S. Solecki, {\em $G_\delta$ ideals of compact sets}, J. Eur. Math. Soc. 13 (2011), 853--882.

\bibitem{ST1} S. Solecki, S. Todorcevic, {\em Cofinal types of topological directed orders}, Ann. Inst.
Fourier (Grenoble) 54 (2004), 1877--1911.

\bibitem{ST2} S. Solecki, S. Todorcevic, {\em Avoiding families and Tukey functions on the nowhere dense ideal},
J. Inst. Math. Jussieu 10 (2011), 405--435.

\bibitem{T1} S. Todorcevic, {\em A classification of transitive relations on $\omega_1$},
Proc. London Math. Soc. 73 (1996), 501--533.

\bibitem{T2} S. Todorcevic, {\em Analytic gaps}, Fund. Math. 150 (1996), 55--66.

\bibitem{T} S. Todorcevic, {\em Basis problems in combinatorial set theory}, in {\em Proceedings of the
International Congress of Mathematicians}, Vol. II, Doc. Math. 1998, Extra Vol. II, 43--52.

\bibitem{Tu} J.W. Tukey, {\em Convergence and Uniformity in Topology}, Annals of Mathematics Studies, 2, Princeton University Press, 1940.
\end{thebibliography}
\end{document}